\documentclass{amsart}

\usepackage[utf8]{inputenc} 
\usepackage[english]{babel} 
\usepackage{mathrsfs}

\theoremstyle{plain} 
\newtheorem{thm}{Theorem} 
\newtheorem{lem}[thm]{Lemma} 
\newtheorem{cor}[thm]{Corollary}
\newtheorem*{lem*}{Lemma}

\theoremstyle{definition} 
\newtheorem*{defn}{Definition}

\theoremstyle{remark} 
\newtheorem*{rem}{Remark}

\begin{document} 

\title{Composition Operators on Bohr--Bergman Spaces of Dirichlet Series} 
\date{\today} 

\author{Maxime Bailleul} 
\address{Univ Lille-Nord-de-France UArtois, Laboratoire de Math\'ematiques de Lens EA~2462, F\'ed\'eration CNRS Nord-Pas-de-Calais FR~2956, F-62\kern 1mm 300 Lens, France} 
\email{maxime.bailleul@euler.univ-artois.fr} 

\author{Ole Fredrik Brevig} 
\address{Department of Mathematical Sciences, Norwegian University of Science and Technology (NTNU), NO-7491 Trondheim, Norway} 
\email{ole.brevig@math.ntnu.no}

\thanks{The second named author is supported by Grant 227768 of the Research Council of Norway.}

\begin{abstract}
	For $\alpha \in \mathbb{R}$, let $\mathscr{D}_\alpha$ denote the scale of Hilbert spaces consisting of Dirichlet series $f(s) = \sum_{n=1}^\infty a_n n^{-s}$ that satisfy $\sum_{n=1}^\infty |a_n|^2/[d(n)]^\alpha < \infty$. The Gordon--Hedenmalm Theorem on composition operators for $\mathscr{H}^2=\mathscr{D}_0$ is extended to the Bergman case $\alpha>0$. These composition operators are generated by functions of the form $\Phi(s) = c_0 s + \varphi(s)$, where $c_0$ is a nonnegative integer and $\varphi(s)$ is a Dirichlet series with certain convergence and mapping properties. For the operators with $c_0=0$ a new phenomenon is discovered: If $0 < \alpha < 1$, the space $\mathscr{D}_\alpha$ is mapped by the composition operator into a smaller space in the same scale. When $\alpha > 1$, the space $\mathscr{D}_\alpha$ is mapped into a larger space in the same scale. Moreover, a partial description of the composition operators on the Dirichlet--Bergman spaces $\mathscr{A}^p$ for $1 \leq p < \infty$ are obtained, in addition to new partial results for composition operators on the Dirichlet--Hardy spaces $\mathscr{H}^p$ when $p$ is an odd integer.
\end{abstract}


\maketitle

% INTRODUCTION
\section{Introduction} \label{sec:intro} 
A theorem of Gordon and Hedenmalm \cite{gordonhedenmalm} describes the composition operators on the Hardy space $\mathscr{H}^2$ of ordinary Dirichlet series with square summable coefficients. In the present work, we consider a scale of weighted Hilbert spaces of Dirichlet series that are analogues to the weighted Bergman spaces in the unit disc, and extend the Gordon--Hedenmalm Theorem to these spaces. To obtain this result, we will rely in part on the tools from \cite{gordonhedenmalm}, but also on new techniques where we use certain averages of twisted Dirichlet series and twisted composition operators.

We let $\mathbb{C}_\theta$ denote the half-plane of complex numbers $s = \sigma+it$ with $\sigma>\theta$. The Dirichlet series in $\mathscr{H}^2$ represent analytic functions in $\mathbb{C}_{1/2}$. A slight strengthening of the Gordon--Hedenmalm Theorem \cite{queffelecseip} states that $\Phi \colon \mathbb{C}_{1/2} \to \mathbb{C}_{1/2}$ generates a composition operator on $\mathscr{H}^2$ if and only if $\Phi$ is a member of the following class: 
\begin{defn}
	The \emph{Gordon--Hedenmalm class}, denoted $\mathscr{G}$, is the set of functions $\Phi\colon\mathbb{C}_{1/2}\to\mathbb{C}_{1/2}$ of the form 
	\begin{equation}
		\label{eq:compclass} \Phi(s) = c_0 s + \sum_{n=1}^\infty c_n n^{-s} = c_0s + \varphi(s), 
	\end{equation}
	where $c_0$ is a nonnegative integer called the \emph{characteristic} of $\Phi$. The Dirichlet series $\varphi$ converges uniformly in $\mathbb{C}_\varepsilon$ $(\varepsilon>0)$ and has the following mapping properties:
	
	{\normalfont(a)} If $c_0=0$, then $\varphi(\mathbb{C}_0)\subset \mathbb{C}_{1/2}$.
	
	{\normalfont(b)} If $c_0\geq1$, then either $\varphi\equiv 0$ or $\varphi(\mathbb{C}_0)\subset \mathbb{C}_0$. 
\end{defn}

It is important to notice that even though $\Phi \in \mathscr{G}$ is a mapping from $\mathbb{C}_{1/2}$ to itself, $\Phi$ is always defined on all of $\mathbb{C}_0$. By the \emph{mapping properties} of $\Phi \in \mathscr{G}$, we mean both the fact that $\Phi(\mathbb{C}_{1/2})\subseteq \mathbb{C}_{1/2}$ and condition (a) or (b) above, depending on the characteristic.

Now, let $d(n)$ denote the number of divisors of the integer $n$. For real numbers $\alpha$, we consider the following scale of Hilbert spaces of Dirichlet series:
\begin{equation}
	\mathscr{D}_\alpha = \left\{f(s) = \sum_{n=1}^\infty a_n n^{-s}\,:\, \|f\|_{\mathscr{D}_\alpha} = \left(\sum_{n=1}^\infty \frac{|a_n|^2}{[d(n)]^\alpha}\right)^\frac{1}{2}<\infty\right\}.  \label{eq:Dalpha} 
\end{equation}
The Hardy space $\mathscr{H}^2$ is the case $\alpha=0$. We will let $\mathscr{A}^2$ denote the case $\alpha=1$, which is a Dirichlet series analogue to the classical unweighted Bergman space of the unit disc, $A^2(\mathbb{D})$. Since $d(n)=\mathcal{O}(n^\varepsilon)$ for every $\varepsilon>0$, see \cite[Thm. 315]{hardywright}, the Cauchy--Schwarz inequality implies that every $f\in\mathscr{D}_\alpha$ is analytic in $\mathbb{C}_{1/2}$. Our main result is:

\begin{thm}
	\label{thm:compo} Let $\alpha>0$. A function $\Phi \colon \mathbb{C}_{1/2} \to \mathbb{C}_{1/2}$ generates a composition operator $\mathscr{C}_\Phi \colon \mathscr{D}_\alpha \to \mathscr{D}_\beta$, where
	\[\beta = \begin{cases}
		2^\alpha-1 & \text{ if } c_0 = 0 \\
		\alpha & \text{ if } c_0 \geq 1
	\end{cases},\]
	 if and only if $\Phi\in\mathscr{G}$. Moreover, if $c_0 \geq 1$ the operator is a contraction.
\end{thm}

The fact that $\beta=\alpha$ is optimal when $c_0\geq1$ is obvious, by considering the composition operator generated by $\Phi(s) = c_0 s$, and the Dirichlet series
\[f(s) = \sum_{k=0}^\infty (k+1)^\frac{\alpha-1}{2}\, 2^{-ks}.\]
Whether $\beta=2^\alpha-1$ is optimal when $c_0=0$ is not clear. However, when $0 < \alpha < 1$ we note that $\mathscr{D}_\alpha$ is mapped by $\mathscr{C}_\Phi$ into the smaller space $\mathscr{D}_\beta$. When $\alpha>1$, it is mapped into the larger space $\mathscr{D}_\beta$. The only cases where $2^\alpha-1=\alpha$ is when $\alpha=0$ or $\alpha=1$, which corresponds to the spaces $\mathscr{H}^2$ and $\mathscr{A}^2$, respectively.

\begin{rem}
	The scale of spaces $\mathscr{D}_\alpha$ is a Dirichlet series version of the classical Dirichlet scale of Hilbert spaces in the unit disc, $D_\alpha(\mathbb{D})$, as defined in \eqref{eq:Dscale}. Every analytic function $\psi \colon \mathbb{D}\to\mathbb{D}$ generates a composition operator on $D_\alpha(\mathbb{D})$ when $\alpha>0$, see \cite[Ch.~11]{zhu}. The functions $\Phi \in \mathscr{G}$ with characteristic $0$ do not fix $+\infty$, and are thus analogues to the functions $\psi$ which do not fix the origin. Interestingly, the phenomenon discovered above does not appear in the classical situation, where such composition operators generally map $D_\alpha(\mathbb{D})$ to $D_\alpha(\mathbb{D})$.
\end{rem}

The remainder of this paper is organized as follows: 
\begin{itemize}
	\item Section~\ref{sec:prelim} includes some preliminary results pertaining to $\mathscr{D}_\alpha$ and on the functions in the class $\mathscr{G}$.
	\item Section~\ref{sec:proof} is devoted to the proof of Theorem~\ref{thm:compo}.
	\item Section~\ref{sec:rem} consists of a partial description of the composition operators on $\mathscr{A}^p$, for $1 \leq p <\infty$, in addition to new partial results regarding composition operators on $\mathscr{H}^p$, when $p$ is an odd integer.
\end{itemize}

We will use the notation $f(x) \ll g(x)$ when there is some constant $C>0$ such that $|f(x)| \leq C |g(x)|$. If both $f(x) \ll g(x)$ and $g(x) \ll f(x)$ hold, we will write $f(x) \asymp g(x)$.

% PRELIMINARIES
\section{Preliminaries} \label{sec:prelim}
\subsection{} One of the most important tools in the study of function spaces of Dirichlet series is the \emph{Bohr lift}, which was introduced by H. Bohr \cite{bohr}. We will lift the Dirichlet series
\begin{equation} \label{eq:diriseri} 
	f(s) = \sum_{n=1}^\infty a_n n^{-s}, 
\end{equation}
to a function on the polydisc $\mathbb{D}^\infty = \{z = (z_1,\,z_2,\,z_3,\,\ldots\,) \, : \, |z_j|<1\}$. Writing $n$ as a product of its prime factors,
\begin{equation} \label{eq:primefact} 
	n = \prod_{j} p_j^{\kappa_j}, 
\end{equation}
we associate $n$ to the finite multi-index $\kappa(n)=(\kappa_1,\,\kappa_2,\,\kappa_3,\,\ldots\,)$. We will sometimes write $n=p^\kappa$ with \eqref{eq:primefact} in mind. The Bohr lift of $f$, denoted by $\mathscr{B}f$, is the power series
\[(\mathscr{B}f)(z) = \sum_{n=1}^\infty a_n z^{\kappa(n)},\]
in view of \eqref{eq:diriseri} and \eqref{eq:primefact}. The polytorus $\mathbb{T}^\infty$ is the distinguished boundary of $\mathbb{D}^\infty$. As shown in \cite{HLS}, the space $\mathscr{H}^2$ is identified with the Hardy space $H^2(\mathbb{T}^\infty)$ under the Bohr lift. A similar identification is obtained for the spaces $\mathscr{H}^p$ by Bayart \cite{bayart}.

In the introduction, we claimed that $\mathscr{A}^2$ is a natural Dirichlet series analogue to the classical Bergman space in the unit disc, $A^2(\mathbb{D})$. To explain this, we let $F$ be an analytic function in $\mathbb{D}$ with Taylor expansion 
\begin{equation} \label{eq:taylorseries} 
	F(z)=\sum_{k=0}^\infty b_k z^k, 
\end{equation}
and let $m$ denote the normalized Lebesgue measure on $\mathbb{D}$. A standard computation shows that 
\begin{equation} \label{eq:Bergmandisc} 
	\|F\|_{A^2(\mathbb{D})} = \left(\int_{\mathbb{D}} |F(w)|^2\, dm(w)\right)^\frac{1}{2} = \left(\sum_{k=0}^\infty \frac{|b_k|^2}{k+1}\right)^\frac{1}{2}.
\end{equation}
In view of the Bohr lift, we let $d\nu(z) = dm(z_1) \times dm(z_2) \times dm(z_3) \times \cdots$. By applying \eqref{eq:Bergmandisc} in each variable, we compute 
\begin{equation} \label{eq:A2norm} 
	\|f\|_{\mathscr{A}^2} = \left(\int_{\mathbb{D}^\infty} |(\mathscr{B}f)(z)|^2\,d\nu(z)\right)^\frac{1}{2}=\left(\sum_{n=1}^\infty\frac{|a_n|^2}{d(n)}\right)^\frac{1}{2}, 
\end{equation}
since $d(n)= (\kappa_1+1)(\kappa_2+1)(\kappa_3+1)\cdots$ when $n=p^\kappa$. Since its norm is defined through the Bohr lift, we refer to $\mathscr{A}^2$ as a \emph{Bohr--Bergman} space. 

We will consider the \emph{Kronecker flow} of the point $z=(z_1,\,z_2,\,z_3,\,\ldots\,)\in \mathbb{C}^\infty$, which is given by
\[\mathscr{T}_t(z)=\left(2^{-it}z_1,\,3^{-it}z_2,\,5^{-it}z_3,\,\ldots\,\right),\qquad t \in \mathbb{R}.\]
The Kronecker flow is simply a rotation in each variable, and defines an ergodic flow on $\mathbb{T}^\infty$ by Kronecker's Theorem \cite[Ch.~13]{hardywright}.

\subsection{} Let us now turn to the spaces $\mathscr{D}_\alpha$ $(\alpha>0)$ as defined in \eqref{eq:Dalpha}. They are Dirichlet series analogues to the spaces $D_\alpha(\mathbb{D})$. The latter space contains the functions $F$ that are analytic in $\mathbb{D}$ and satisfy $\|F\|_{D_\alpha(\mathbb{D})}<\infty$, where
\begin{equation} \label{eq:Dscale}
	\|F\|_{D_\alpha(\mathbb{D})} = \left(\int_{\mathbb{D}} |F(w)|^2\alpha\left(1-|w|^2\right)^{\alpha-1}\,dm(w)\right)^\frac{1}{2} \asymp \left(\sum_{k=1}^\infty \frac{|a_k|^2}{(k+1)^\alpha}\right)^\frac{1}{2}.
\end{equation}
In $\mathbb{D}$, we can use equivalent norms and get the same set of functions. However, since we will use an infinite number of variables, exact equality is needed. Thus, if we were to define a space of Dirichlet series using the infinite product of the Möbius invariant measure $\alpha\left(1-|w|^2\right)^{\alpha-1}\,dm(w)$ in a similar manner to the definition of \eqref{eq:A2norm}, this \emph{would not} be the space $\mathscr{D}_\alpha$.

\begin{lem} \label{lem:Dalpha} 
	Let $\alpha>0$. There exists a rotationally invariant probability measure $\nu_\alpha$ on $\mathbb{D}^\infty$ such that for $f$ of the form \eqref{eq:diriseri} we have
	\[\int_{\mathbb{D}^\infty} |(\mathscr{B}f)(z)|^2\,d\nu_\alpha(z) = \sum_{n=1}^\infty \frac{|a_n|^2}{[d(n)]^\alpha}.\]
	Moreover, for any probability measure $\lambda$ on $\mathbb{R}$ we have 
	\begin{equation} \label{eq:Dalphanorm} 
		\|f\|_{\mathscr{D}_\alpha}^2 = \int_{\mathbb{D}^\infty} \int_\mathbb{R} |(\mathscr{B}f)(\mathscr{T}_t \chi)|^2\,d\lambda(t)d\nu_\alpha(\chi).
	\end{equation}
	\begin{proof}
		For $F$ of the form \eqref{eq:taylorseries} we want to find a radial probability measure $m_\alpha$ such that 
		\begin{equation} \label{eq:weightdemand} 
			\int_{\mathbb{D}} |F(w)|^2 \, dm_\alpha(w) = \sum_{k=0}^\infty \frac{|b_k|^2}{(k+1)^\alpha}.
		\end{equation}
		To this end, we introduce
		\[dm_\alpha(w) = \left(\log{\frac{1}{|w|^2}}\right)^{\alpha-1} \, \frac{dm(w)}{\Gamma(\alpha)},\]
		and define $d\nu_\alpha(z) = dm_\alpha(z_1)\times dm_\alpha(z_2)\times dm_\alpha(z_3)\times \cdots$. We then use $d\nu_\alpha(z)$ to define $\mathscr{D}_\alpha$ in a similar way to \eqref{eq:A2norm}. The measure $d\nu_\alpha(z)$ is clearly rotationally invariant. Hence \eqref{eq:Dalphanorm} follows by Fubini's theorem and the fact that the Kronecker flow is a rotation in each variable. 
	\end{proof}
\end{lem}

\subsection{} Let us consider the following multiplier problem: Let $\mathscr{M}(\mathscr{D}_\alpha)$ denote the collection of analytic functions $m \colon \mathbb{C}_{1/2} \to \mathbb{C}$ so that $mf \in \mathscr{D}_\alpha$ for every $f\in\mathscr{D}_\alpha$. In \cite{HLS}, the space $H^\infty(\mathbb{D}^\infty)$ is shown to be isometrically isomorphic to the space
\[\mathscr{H}^\infty = \left\{f(s) = \sum_{n=1}^\infty a_n n^{-s}\,:\, \|f\|_\infty = \sup_{\sigma>0}|f(s)|<\infty\right\},\]
through the Bohr lift. In particular, this easily implies that $\mathscr{H}^\infty\subseteq\mathscr{M}(\mathscr{D}_\alpha)$, and that $\|m\| \leq \|m\|_\infty$, where the former denotes the norm of $m$ as a multiplier. The following result has previously been observed by Olsen \cite{olsenlocal}, but we include a short proof. A similar result is also proved in \cite[Thm.~11.21]{bailleul}. A general result on multipliers of weighted Hilbert spaces of Dirichlet series with multiplicative weights was recently obtained by Stetler \cite{stetler}.

\begin{thm} \label{thm:multiplier}
		Let $\alpha>0$. The multiplier algebra $\mathscr{M}(\mathscr{D}_\alpha)$ is $\mathscr{H}^\infty$.
		\begin{proof}
			We have already observed that $\mathscr{H}^\infty \subseteq \mathscr{M}(\mathscr{D}_\alpha)$, and that $\|m\|\leq \|m\|_\infty$. For the other inclusion, we observe that since $1 \in \mathscr{D}_\alpha$ we have $m\in \mathscr{D}_\alpha$, and hence $m$ is a Dirichlet series. In fact, $m^j$ is in $\mathscr{D}_\alpha$ for every $j \in \mathbb{N}$, and $\|m^j\|_{\mathscr{D}_\alpha}\leq\|m\|^j$. This implies that
			\[\|m\|_\infty = \sup_{z \in \mathbb{D}^\infty} |(\mathscr{B}m)(z)| = \lim_{j \to \infty} \left(\int_{\mathbb{D}^\infty} |(\mathscr{B}m)(z)|^{2j} d\nu_\alpha(z)\right)^\frac{1}{2j} \leq \|m\|,\]
			which concludes the proof.
		\end{proof}
\end{thm}

Curiously, even though we initially only require multipliers to be defined in $\mathbb{C}_{1/2}$, they are automatically defined in the larger half-plane $\mathbb{C}_0$. This phenomenon is typical for function spaces of Dirichlet series, and it also appears in the definition of the Gordon--Hedenmalm class $\mathscr{G}$. Furthermore, the phenomenon indicates that we need to obtain properties of $\mathscr{D}_\alpha$ in both half-planes $\mathbb{C}_0$ and $\mathbb{C}_{1/2}$. 

\subsection{} First, for $\chi=(\chi_1,\,\chi_2,\,\chi_3,\,\ldots\,)\in\mathbb{C}^\infty$ we define a completely multiplicative function by $\chi(n)=\chi^\kappa$, when $n=p^\kappa$. For $f$ of the form \eqref{eq:diriseri}, we consider the \emph{twisted Dirichlet series}
\begin{equation}  \label{eq:diritwist} 
	f_\chi(s) = \sum_{n=1}^\infty a_n \chi(n) n^{-s}. 
\end{equation}
If $\chi\in\mathbb{T}^\infty$, then $f_\chi$ is just a \emph{vertical limit function} of $f$, as considered in \cite{gordonhedenmalm,HLS}. We will also consider the \emph{shifted Dirichlet series} $f_\delta(s)=f(s+\delta)$. If a Dirichlet series is both twisted and shifted, we will write $f_{\delta,\chi}$. Observe that for those $\chi\in\mathbb{C}^\infty$ and $s=\sigma+it$ where the series \eqref{eq:diritwist} converges we have
\begin{equation} \label{eq:wehaveeq}
	f_\chi(s) = (\mathscr{B}f_\sigma)(\mathscr{T}_t\chi).
\end{equation}
In particular, \eqref{eq:wehaveeq} is true for $s \in \mathbb{C}_{1/2}$ and $\chi \in \overline{\mathbb{D}^\infty}$. This follows by the Cauchy--Schwarz inequality and the fact that $|\chi(n)|\leq1$, since $\chi \in \overline{\mathbb{D}^\infty}$. To further extend the validity of \eqref{eq:wehaveeq}, we will use the following result:

\begin{lem*}[Rademacher--Menchov]
	Let $(X,\mu)$ be a probability space and suppose that $\{e_n\}_{n=1}^\infty$ is an orthonormal sequence in $L^2(X)$. If $\sum_{n=1}^\infty |c_n|^2(\log{n})^2 < \infty$ then the series 
	\[\sum_{n=1}^\infty c_n e_n(x)\]
	converges for $\mu$-almost every $x \in X$.
	\begin{proof}
		A proof may be found in \cite[p. 42]{olevskii}
	\end{proof}
\end{lem*}

\begin{lem} \label{lem:twistext}
	Let $f \in \mathscr{D}_\alpha$. The Dirichlet series $f_\chi$ as defined by \eqref{eq:diritwist} converges in $\mathbb{C}_0$ for almost every $\chi\in\mathbb{D}^\infty$, with respect to $\nu_\alpha$. 
	\begin{proof}
		We shall use the the Rademacher--Menchov Lemma on the space $L^2(\mathbb{D}^\infty,\nu_\alpha)$. Hence we let 
		\[e_n(\chi)=[d(n)]^{\alpha/2}\chi(n)\qquad \text{and} \qquad c_n =  a_n [d(n)]^{-\alpha/2}n^{-s},\] 
		for $f$ of the form \eqref{eq:diriseri}. When $\sigma>0$, the Rademacher--Menchov Lemma implies that $f_\chi(s)$ converges for almost every $\chi\in\mathbb{D}^\infty$. It is well known that if a Dirichlet series converges in a point $s_0 = \sigma_0 + i t_0$, then it also converges for every $s$ with $\sigma>\sigma_0$. In particular, we may take $s = 1/j$ and conclude that for almost every $\chi\in\mathbb{D}^\infty$, the series $f_\chi$ converges in $\mathbb{C}_{1/j}$. Clearly this implies that $f_\chi$ converges in $\mathbb{C}_0$ for almost every $\chi\in\mathbb{D}^\infty$, since the union of sets of zero measure has zero measure. 
	\end{proof}
\end{lem}

Thus we conclude that \eqref{eq:wehaveeq} is true for almost every $\chi \in \mathbb{D}^\infty$ and every $s \in \mathbb{C}_0$. Now, we would like a version of \eqref{eq:Dalphanorm} for $f_\chi(s)$, but in general we do not know anything about the existence of $f_\chi(it)$. However, in light of Lemma~\ref{lem:twistext} we may combine \eqref{eq:Dalphanorm} with \eqref{eq:wehaveeq} and use Fubini's theorem to obtain the formula
\begin{equation} \label{eq:twistnorm}
	\|f\|_{\mathscr{D}_\alpha} = \lim_{\sigma\to 0^+}\left(\int_{\mathbb{D}^\infty} \int_{\mathbb{R}} |f_\chi(\sigma+it)|^2\,d\lambda(t)d\nu_\alpha(\chi)\right)^\frac{1}{2},
\end{equation}
which will be sufficient for our purposes.

\subsection{} Let us turn to the half-plane $\mathbb{C}_{1/2}$. Any function in $\mathscr{D}_\alpha$ can be expressed as a limit of Dirichlet polynomials, with convergence in the norm of $\mathscr{D}_\alpha$. Now, let $\tau \colon \mathbb{C}_{1/2}\to\mathbb{D}$ be the conformal mapping defined by
\begin{equation} \label{eq:taudef}
	\tau(s) = \frac{s-3/2}{s+1/2}.
\end{equation}
For $\beta>0$ we let $D_{\beta,\mathrm{i}}(\mathbb{C}_{1/2})$ denote the space defined by the pull-back of $\tau$ from $D_\beta(\mathbb{D})$. This means that $D_{\beta,\mathrm{i}}(\mathbb{C}_{1/2})$ consists of functions $f$ that are analytic in $\mathbb{C}_{1/2}$ and finite with respect to the norm
\[\|f\|_{D_{\beta,\mathrm{i}}(\mathbb{C}_{1/2})}=\|f\circ \tau^{-1}\|_{D_\beta(\mathbb{D})} = 4^\beta\beta \int_{\mathbb{C}_{1/2}} |f(s)|^2 \left(\sigma-\frac{1}{2}\right)^{\beta-1} \frac{dm(s)}{|s+1/2|^{2\beta+2}}.\]
We shall need the following embedding result:
\begin{lem} \label{lem:localembed} 
	$\mathscr{D}_\alpha$ is continuously embedded into $D_{\beta,\mathrm{i}}(\mathbb{C}_{1/2})$, where $\beta=2^\alpha-1$.
	\begin{proof}
		The embedding can be deduced from the corresponding local embedding (see Theorem~1 and Example~4 in \cite{olsenlocal}) by straightforward estimates.
	\end{proof}
\end{lem}

\begin{rem}
	The reason behind the relationship $\beta=2^\alpha-1$ is the classical asymptotic formula
	\[\sum_{n\leq x} [d(n)]^\alpha \asymp x(\log{x})^\beta,\]
	which is due to Ramanujan \cite{ramanujan} and Wilson \cite{wilson}. The embedding is optimal, in the sense that we cannot replace $\beta=2^\alpha-1$ with any smaller value \cite{olsenlocal}. 
\end{rem}

\subsection{} We shall now consider twisted composition operators by extending \eqref{eq:diritwist} to functions $\Phi\in\mathscr{G}$ by defining 
\[\Phi_\chi(s) = c_0s + \varphi_\chi(s).\]
We will also write $\Phi_{\chi,\delta}(s)=\Phi_\chi(s+\delta)$. We expect that the functions $\Phi_\chi$ behave similarly to $\Phi$, and the following result shows that the mapping and convergence properties of $\Phi$ are retained in $\Phi_\chi$. This is crucial, as we need to replace the ``limit measure'' used in \cite{gordonhedenmalm} with certain mean values of the composition operators generated by $\Phi_\chi$.

\begin{lem} \label{lem:twistmap} 
	Suppose that $\Phi \in \mathscr{G}$. Then $\Phi_\chi \in \mathscr{G}$ for any $\chi \in \overline{\mathbb{D}^\infty}$. 
	\begin{proof}
		As explained in \cite{gordonhedenmalm} (see Proposition 4.1), if $\chi \in \mathbb{T}^\infty$ the function $\varphi_\chi$ is a normal limit of vertical translates of $\varphi$. Since vertical translation does not change the mapping properties of $\varphi$, neither does the application of $\chi \in \mathbb{T}^\infty$ to $\Phi$. To extend this into $\mathbb{D}^\infty$, we let $\chi \in \mathbb{T}^\infty$ and consider
		\[\chi(w) = (w,\,\chi_2,\,\chi_3,\,\chi_4,\,\ldots\,).\]
		We fix $s \in \mathbb{C}_0$ and apply the maximum modulus principle on $\mathbb{D}$ to the one-variable analytic function $F(w) = \exp\left(-\varphi_{\chi(w)}(s)\right)$. This implies that $\varphi_{\chi(w)}$ maps $\mathbb{C}_0$ to $\mathbb{C}_\theta$ for any $w \in \mathbb{D}$. We employ the same procedure to every coordinate of $\chi$.
		
		By a theorem of Bohr \cite{bohruni}, the Dirichlet series $\varphi$ converges uniformly in any closed half-plane where it can be represented as a bounded analytic function. Since vertical translation of $\varphi$ does not change boundedness or analyticity, it is clear that $\varphi_\chi$ converges uniformly where $\varphi$ converges uniformly, when $\chi\in \mathbb{T}^\infty$. By a similar maximum modulus argument as above with $G(w)= \varphi_{\chi(w)}(s)$, this can be extended into $\mathbb{D}^\infty$. 
	\end{proof}
\end{lem}

The following result is a version of \cite[Prop.~4.3]{gordonhedenmalm} for Dirichlet polynomials. Our version plays a prominent role in the proof of the sufficiency part of Theorem~\ref{thm:compo}, whereas in \cite{gordonhedenmalm} the corresponding result is only used in the proof of the necessity part.
\begin{lem} \label{lem:twistcomp} 
	Suppose that $\Phi\in\mathscr{G}$. For every Dirichlet polynomial $f$, every $\chi \in \overline{\mathbb{D}^\infty}$ and every $s \in \mathbb{C}_0$, we have 
	\begin{equation} \label{eq:twistcomp}
		\left(f\circ \Phi\right)_\chi(s) = \left(f_{\chi^{c_0}}\circ \Phi_{\chi}\right)(s)
	\end{equation}
	where $\chi^{c_0}=\left(\chi_1^{c_0},\,\chi_2^{c_0},\,\chi_3^{c_0},\,\ldots\,\right)$.
	\begin{proof}
		A formal computation extracted from \cite[Sec.4]{gordonhedenmalm} shows that both sides of \eqref{eq:twistcomp} are Dirichlet series, and that they are equal. By the assumptions on $f$ and $\Phi$, it is evident that both sides converge absolutely in (at least) $\mathbb{C}_1$, so \eqref{eq:twistcomp} is valid there. The right hand side represents a bounded analytic function in $\mathbb{C}_0$, so by Bohr's theorem (see the proof of Lemma~\ref{lem:twistmap}) and the identity principle this extends to $\mathbb{C}_0$.
	\end{proof}
\end{lem}

\begin{rem}
	It is possible to extend Proposition~4.3 in \cite{gordonhedenmalm} to $\mathscr{D}_\alpha$, either by a variation of the argument given above or by the argument used in \cite{gordonhedenmalm}. In the latter case, we appeal to the maximum modulus principle when passing from $\chi \in \mathbb{T}^\infty$ to $\chi \in \overline{\mathbb{D}^\infty}$.
\end{rem}

% PROOF
\section{Proof of Theorem~\ref{thm:compo}} \label{sec:proof} 
\noindent The proof of Theorem~\ref{thm:compo} can essentially be split into three distinct parts. The first two parts are easy to obtain from \cite{gordonhedenmalm}, while our new techniques will be needed in the third. Note that these new techniques can also be applied to prove the corresponding part of the Gordon--Hedenmalm Theorem for $\mathscr{H}^2$.

\subsection{} The first part is the so-called ``arithmetical condition'', which demands that $\Phi$ is of the form \eqref{eq:compclass} to ensure that $f\circ\Phi$ is a somewhere a convergent Dirichlet series. The proof for $\mathscr{D}_\alpha$ translates directly from the work of Gordon and Hedenmalm on $\mathscr{H}^2$, see \cite[Thm. A]{gordonhedenmalm}.

\subsection{} The second part is the necessity of the mapping and convergence properties of $\Phi$. The argument given for $\mathscr{H}^2$ in \cite{gordonhedenmalm} is quite general, and applies almost directly to $\mathscr{D}_\alpha$. We need only observe that Lemma~\ref{lem:twistext} still holds with $\mathbb{D}^\infty$ replaced by $\mathbb{T}^\infty$ and $\nu_\alpha$ replaced by the Haar measure $\mu$ of $\mathbb{T}^\infty$. The argument given in \cite{gordonhedenmalm} will then apply line for line provided we can prove the following lemma:
\begin{lem} \label{lem:primefunc} 
	There is a function $f \in \mathscr{D}_\alpha$ with the following properties:
	
	{\normalfont (i)} For almost every $\chi \in \mathbb{T}^\infty$, $f_\chi$ converges in $\mathbb{C}_0$ and cannot be analytically continued to any larger domain.
	
	{\normalfont (ii)} For at least one $\chi \in \mathbb{T}^\infty$, $f_\chi$ converges in $\mathbb{C}_{1/2}$ and cannot be analytically continued to any larger domain.
	
	\begin{proof}
		The function in question is
		\begin{align*}
			f(s) &= \sum_p c(p) p^{-s}, \\
			c(x) &= \frac{1}{\sqrt{x}\log{x}},
		\end{align*}
		where the sum is taken over the prime numbers. As verified in \cite[Lem.~9]{bayart}, $f$ satisfies the required properties. It is also evident that $f\in \mathscr{D}_\alpha$, since $d(p)=2$.
	\end{proof}
\end{lem}
Moreover, the strengthening of the necessity argument due to Queff{\'e}lec and Seip also applies to $\mathscr{D}_\alpha$, see Section 3 and in particular Theorem~3.1 in \cite{queffelecseip}.

\subsection{} The third and final part of the proof is the sufficiency of the mapping properties. This is where the ``limit measure procedure'' of \cite{gordonhedenmalm} does not apply, and we have to find new techniques. Curiously, the cases $c_0=0$ and $c_0 \geq 1$ are handled quite differently: The first case is done by integration over $\mathbb{D}$ and $\mathbb{T}^\infty$, and the second case is done by integration over $\mathbb{T}$ and $\mathbb{D}^\infty$.

We will need several versions of \emph{Littlewood's subordination principle}, which in its most basic form can be stated as follows: If $\psi \colon \mathbb{D}\to\mathbb{D}$ is analytic with $\psi(0)=0$, then
\[\int_0^{2\pi} |F(\psi(re^{i\theta}))|^p\,d\theta \leq \int_0^{2\pi} |F(re^{i\theta})|^p\,d\theta\]
for every $F$ that is analytic in $\mathbb{D}$, every $0<p<\infty$ and every $0<r<1$. There are various versions of this principle for the various function spaces in $\mathbb{D}$, and we refer generally to \cite[Ch.~11]{zhu}.

\begin{proof}[Proof of sufficiency when $c_0=0$]
	Let $\Phi \in \mathscr{G}$ with $c_0=0$, and consider the map $\tau \colon \mathbb{C}_{1/2} \to \mathbb{D}$ as defined in \eqref{eq:taudef}. We will use the following version of Littlewood's principle: Let $\psi \colon \mathbb{D} \to \mathbb{C}_{1/2}$ be analytic. For every $\beta>0$, the function $\psi$ generates a composition operator $\mathscr{C}_\psi$ from $D_{\beta,\mathrm{i}}(\mathbb{C}_{1/2})$ to $D_\beta(\mathbb{D})$, and
	\[\|\mathscr{C}_\psi\| \leq \left(\frac{1+|\tau(\psi(0))|}{1-|\tau(\psi(0))|}\right)^{(1+\beta)/2}.\]
	Fix some $s \in \mathbb{C}_0$ and some $\chi \in \mathbb{T}^\infty$, and define $\psi_\chi \colon \mathbb{D} \to \mathbb{C}_{1/2}$ by
	\begin{equation} \label{eq:psidef}
		\psi_\chi(w) = \Phi_{w\chi}(s) = \sum_{n=1}^\infty c_n w^{\Omega(n)} \chi(n) n^{-s},
	\end{equation}
	in view of \eqref{eq:compclass}. Here $\Omega(n)$ denotes the number of prime factors of the integer $n$. In \eqref{eq:psidef}, we have introduced the notation $w \chi = (w\chi_1,\,w\chi_2,\,w\chi_3,\,\ldots\,)$. Now, let $f$ be a Dirichlet polynomial, and since $f \circ\Phi$ is a Dirichlet series by Lemma~\ref{lem:twistcomp}, we may write
	\[f\circ\Phi(s) = \sum_{n=1}^\infty b_n n^{-s}.\]
	By appealing again to Lemma~\ref{lem:twistcomp}, we compute
	\begin{equation} \label{eq:psicomp}
		f\circ\psi_\chi(w) = f \circ \Phi_{w\chi}(s) = \left(f\circ \Phi\right)_{w\chi}(s) = \sum_{n=1}^\infty b_n w^{\Omega(n)} \chi(n) n^{-s}.
	\end{equation}
	Clearly $\psi_\chi(0) = c_1$ for every $s \in \mathbb{C}_0$ and every $\chi \in \mathbb{T}^\infty$, and hence Littlewood's principle implies that
	\begin{equation} \label{eq:littlec0}
		\|f\circ\psi\|_{D_\beta(\mathbb{D})}^2 \leq \left(\frac{1+|\tau(c_1)|}{1-|\tau(c_1)|}\right)^{1+\beta}\|f\|_{D_{\beta,\mathrm{i}}(\mathbb{C}_{1/2})}^2.
	\end{equation}
	We observe that there is no $\chi$ on the right hand side of \eqref{eq:littlec0}. We therefore integrate over $\mathbb{T}^\infty$ with respect to the Haar measure $d\mu(\chi)$ on both sides. Clearly, the right hand side does not change. The left hand side of \eqref{eq:littlec0} may be computed using the representation \eqref{eq:psicomp} with Fubini's theorem and \eqref{eq:Dscale}: 
	\[\int_{\mathbb{T}^\infty} \int_{\mathbb{D}} \left|\sum_{n=1}^\infty b_n w^{\Omega(n)} \chi(n) n^{-s}\right|^2\,\beta(1-|w|^2)^{\beta-1}dm(w)d\mu(\chi) \asymp \sum_{n=1}^\infty \frac{|b_n|^2}{[1+\Omega(n)]^\beta}n^{-2\sigma}.\]
	There is no $\sigma$ on the right hand side of \eqref{eq:littlec0}, and hence we may let $\sigma \to 0^+$ on the left hand side to obtain
	\begin{equation} \label{eq:Omegapull} 
		\sum_{n=1}^\infty \frac{|b_n|^2}{[1+\Omega(n)]^\beta} \leq \left(\frac{1+|\tau(c_1)|}{1-|\tau(c_1)|}\right)^{1+\beta}\|f\|_{D_{\beta,\mathrm{i}}(\mathbb{C}_{1/2})}^2. 
	\end{equation}
	The proof is completed by using the fact that $1+ \Omega(n) \leq d(n)$ on the left hand side of \eqref{eq:Omegapull} and Lemma~\ref{lem:localembed} on the right hand side of \eqref{eq:Omegapull}.
\end{proof}

\begin{rem}
	By being more precise, it is possible to obtain $\tau(c_1)=0$ in \eqref{eq:Omegapull}. However, this would not imply that $\mathscr{C}_\Phi$ maps $\mathscr{D}_\alpha$ contractively to $\mathscr{D}_\beta$, since a constant appears when using Lemma~\ref{lem:localembed}. This is as expected, since the point at infinity is not fixed by $\Phi$. Curiously, since $1+\Omega(n)=d(n)$ only when $n=p^k$ for a prime $p$, we still have a certain contractivity. In fact, what we have proved is that $\mathscr{C}_\Phi$ maps $\mathscr{D}_\alpha$ into $\mathscr{H}_\beta$, the Hilbert space of Dirichlet series of the form \eqref{eq:diriseri} that satisfy
	\[\sum_{n=1}^\infty \frac{|a_n|^2}{[1+\Omega(n)]^\beta} < \infty.\]
\end{rem}

\begin{proof}[Proof of sufficiency and contractivity when $c_0 \geq 1$]
	Let $\Phi \in \mathscr{G}$ with $c_0 \geq 1$ and let $\xi>0$ be large, but arbitrary. We will use the following maps from $\mathbb{C}_0$ to $\mathbb{D}$:
	\[\tau_1(s) = \frac{s-c_0\xi}{s+c_0\xi} \qquad \text{ and } \qquad \tau_2(s) = \frac{s - \xi}{s+\xi}.\]
	Let $\chi \in \mathbb{D}^\infty$ and $0 < \sigma \leq 1$ be fixed, and consider $\psi \colon \mathbb{D} \to \mathbb{D}$ defined by
	\[\psi(w) = \left(\tau_1 \circ \Phi_{\chi,\sigma} \circ \tau_2^{-1}\right)(w).\]
	For every $\chi\in\mathbb{D}^\infty$ and every $0 < \sigma \leq 1$, the function $\psi$ is an analytic self-map of $\mathbb{D}$, by Lemma~\ref{lem:twistmap}. Littlewood's subordination principle implies that for every $F \in H^2(\mathbb{D})$ we have
	\begin{equation}  \label{eq:hardysub} 
		\left\|F \circ \psi \right\|_{H^2(\mathbb{D})}^2 \leq \frac{1+|\left(\tau_1 \circ \Phi_{\chi,\sigma} \circ \tau_2^{-1}\right)(0)|}{1-|\left(\tau_1 \circ \Phi_{\chi,\sigma} \circ \tau_2^{-1}\right)(0)|}\|F\|_{H^2(\mathbb{D})}^2.
	\end{equation}
	A direct computation verifies that if $\Phi_{\chi,\sigma}(s) = c_0(s+\sigma) + \varphi_\chi(s+\sigma)$, we have
	\[\left(\tau_1 \circ \Phi_{\chi,\sigma} \circ \tau_2^{-1}\right)(0) = \frac{c_0\sigma+\varphi_\chi(\xi+\sigma)}{c_0(2\xi+\sigma) + \varphi_\chi(\xi+\sigma)}.\]
	Clearly $\varphi_\chi(\xi+\sigma) \to c_1$ as $\xi \to \infty$, and this is uniform in $\chi$ and $\sigma$. Thus, for every $\varepsilon>0$ we may find $\xi$ large, but independent of $\chi$ and $\sigma$, so that 
	\[\frac{1+|\left(\tau_1 \circ \Phi_{\chi,\sigma} \circ \tau_2^{-1}\right)(0)|}{1-|\left(\tau_1 \circ \Phi_{\chi,\sigma} \circ \tau_2^{-1}\right)(0)|} \leq 1 + \varepsilon.\]
	Now, let $f$ be a Dirichlet polynomial, and define
	\[F(w) = \left(f_{\chi^{c_0}}\circ \tau_1^{-1}\right)(w).\]
	Since $f$ is entire and uniformly bounded in $\mathbb{C}_0$, it is clear that $F \in H^\infty(\mathbb{D})$, and hence in $H^2(\mathbb{D})$. Using Lemma~\ref{lem:twistcomp} we obtain
	\[\left(F \circ \psi\right)(w) = \left(f_{\chi^{c_0}} \circ \Phi_{\chi,\sigma} \circ \tau_2^{-1}\right)(w) = \left((f\circ\Phi)_{\chi,\sigma} \circ \tau_2^{-1}\right)(w).\]
	The pull-back of the normalized Lebesgue measure on $\mathbb{T}$ with respect to $\tau_1$ and $\tau_2$ produces the following probability measures on $\mathbb{R}$:
	\[d\lambda_1(t) = \frac{c_0 \xi}{\pi}\,\frac{1}{t^2+(c_0\xi)^2} \qquad \text{ and } \qquad d\lambda_2(t) = \frac{\xi}{\pi}\,\frac{1}{t^2+\xi^2}.\]
	Inserting everything into \eqref{eq:hardysub} we obtain
	\[\int_\mathbb{R} \left|\left(f\circ \Phi\right)_\chi(\sigma+it)\right|^2 d\lambda_1(t) \leq (1+\varepsilon) \int_\mathbb{R} \left|f_{\chi^{c_0}}(it)\right|^2 d\lambda_2(t).\]
	By \eqref{eq:twistnorm}, and by keeping in mind that $f$ is entire and uniformly bounded on $i\mathbb{R}$, we may integrate over $\mathbb{D}^\infty$ with respect to $d\nu_\alpha(\chi)$ and let $\sigma\to0$, to obtain
	\begin{equation} \label{eq:finalest}
		\|f\circ\Phi\|_{\mathscr{D}_\alpha}^2 \leq (1+\varepsilon) \left\|f_{\chi^{c_0}}\right\|_{\mathscr{D}_\alpha}^2 \leq (1+\varepsilon)\left\|f\right\|_{\mathscr{D}_\alpha}^2.
	\end{equation}
	The final inequality in \eqref{eq:finalest} follows from the fact that $d(n^{c_0})\geq d(n)$. Since $\varepsilon>0$ was arbitrary and independent of $f$, the composition operator $\mathscr{C}_\Phi$ is a contraction on $\mathscr{D}_\alpha$.
\end{proof}

% REMARKS
\section{Composition Operators on $\mathscr{A}^p$, $\mathscr{H}^p$ and $\mathcal{A}_\alpha$} \label{sec:rem}
\subsection{} Following the description of the composition operators on $\mathscr{H}^2$ \cite{gordonhedenmalm}, Bayart \cite{bayart} extended the Gordon--Hedenmalm Theorem to $\mathscr{H}^p$, with one exception: The sufficiency of the case (a) is proved only when $p$ is an even integer. Hence the complete description of composition operators on $\mathscr{H}^p$ remains unsolved. 

In view of this, and our results for $\mathscr{D}_\alpha$, it is natural to investigate the composition operators on the Bohr--Bergman spaces $\mathscr{A}^p$, which we for $1 \leq p < \infty$ define as
\[\mathscr{A}^p = \left\{f(s) = \sum_{n=1}^\infty a_n n^{-s}\,:\, \|f\|_{\mathscr{A}^p} = \left(\int_{\mathbb{D}^\infty}|(\mathscr{B}f)(z)|^p\,d\nu(z)\right)^\frac{1}{p}<\infty\right\}.\]
It is convenient to let $\mathscr{A}^\infty=\mathscr{H}^\infty$. Basic properties of the spaces $\mathscr{A}^p$ have been studied in the first named author's thesis \cite[Ch. 11]{bailleul} and in the work of the first named author and Lef{\`e}vre \cite[Sec. 3]{BL2013}.

In particular, we mention that Lemma~\ref{lem:twistext} and hence \eqref{eq:twistnorm} extend to these spaces, mutatis mutandis. Moreover, Lemma~\ref{lem:primefunc} also holds for $\mathscr{A}^p$ by a similar application of Khintchin's inequality for Steinhaus variables as used for $\mathscr{H}^p$ in \cite{bayart}. However, Lemma~\ref{lem:localembed} relies heavily on Hilbert space techniques \cite{olsenlocal}, and the situation for $\mathscr{A}^p$ (and for $\mathscr{H}^p$) is not clear.

With the exception Lemma~\ref{lem:localembed}, our arguments apply almost line for line, and we are able to extend Th{\'e}or{\`e}me 13.6 from \cite{bailleul} and obtain the following result:
\begin{thm} \label{thm:Apcomp}
	Let $1 \leq p < \infty$. Suppose that the function $\Phi \colon \mathbb{C}_{1/2} \to \mathbb{C}_{1/2}$ defines a composition operator $\mathscr{C}_\Phi \colon \mathscr{A}^p \to \mathscr{A}^p$. Then $\Phi \in \mathscr{G}$. Moreover, if $c_0\geq1$ this condition is sufficient and the composition operator is a contraction. If $c_0=0$ and $p=2k$ the condition is sufficient.
	\begin{proof}
		The proof of the ``arithmetical condition'' and the proof of the necessity again follows by the argument in \cite{gordonhedenmalm}. Our proof of the sufficiency when $c_0\geq1$ for $\mathscr{D}_\alpha$ applies line for line, provided we are able to prove the inequality
		\[\int_{\mathbb{D}^\infty} \left|(\mathscr{B}f)(\chi^{c_0})\right|^p\,d\nu(\chi) \leq \int_{\mathbb{D}^\infty} \left|(\mathscr{B}f)(\chi)\right|^p\,d\nu(\chi).\]
		This inequality follows immediately from the fact that the composition operator $F(z) \mapsto F(z^{c_0})$ is contractive on $A^p(\mathbb{D})$ when $c_0\geq 1$. We apply this in each variable with the integral version of Minkowski's inequality. 
		
		To prove the sufficiency $c_0=0$ and $p=2k$, we first observe that for $p=2$, this is just Theorem~\ref{thm:compo} with $\alpha=1$, since $\mathscr{D}_1=\mathscr{A}^2$. This extends immediately to $p=2k$, by the simple fact that for Dirichlet polynomials $f$ we have
		\[\|f\|_{\mathscr{A}^p}^p = \|f^{p}\|_{\mathscr{A}^2}^2,\] 
		and since clearly $(f\circ\Phi)^p(s) = (f^p\circ\Phi)(s)$ for every $s\in\mathbb{C}_0$ and $\Phi \in \mathscr{G}$.
	\end{proof}
\end{thm}

The result of Theorem~\ref{thm:Apcomp} mirrors that of Bayart for $\mathscr{H}^p$, since we are not able to prove sufficiency when $c_0=0$ and $p \neq 2k$. This is not at all surprising, and seems to be due to the fact that we lack local embeddings of $\mathscr{A}^p$ into Bergman spaces in $\mathbb{C}_{1/2}$ when $p\neq2k$, a similar situation to that for $\mathscr{H}^p$. For more on the embedding problem, we refer to \cite{saksmanseip}.

\subsection{} Let us now turn to the spaces $\mathscr{H}^p$, as introduced in \cite{bayart}. We may compute the $\mathscr{H}^p$ norm of a Dirichlet polynomial $f$ in two different ways:
\begin{equation} \label{eq:Hpnorm}
	\|f\|_{\mathscr{H}^p} = \lim_{T\to\infty} \frac{1}{T}\left(\int_0^T |f(it)|^p\,dt\right)^\frac{1}{p} = \left(\int_{\mathbb{T}^\infty} |\mathscr{B}f(z)|^p\,d\mu(z)\right)^\frac{1}{p}.
\end{equation}
Here, $d\mu$ again denotes the normalized Haar measure of $\mathbb{T}^\infty$. The validity of formula \eqref{eq:Hpnorm} follows from the ergodicity of the Kronecker flow on $\mathbb{T}^\infty$. The Kronecker flow is clearly not ergodic on $\mathbb{D}^\infty$ with respect to $d\nu(z)$, so we cannot expect to have a formula of the type \eqref{eq:Hpnorm} for $\mathscr{A}^p$.

Since the situation for composition operators with characteristic $c_0=0$ of $\mathscr{H}^p$ is not clear, we seek a partial result. We want to find a space that is mapped into $\mathscr{H}^p$ by $\mathscr{C}_\Phi$ and a space that $\mathscr{H}^p$ is mapped into by $\mathscr{C}_\Phi$, when $\Phi\in\mathscr{G}$ with $c_0=0$. 

We let $\mathscr{K}=\mathscr{H}^2\odot\mathscr{H}^2$ denote Helson's space. We refer to \cite{helsonbook} and \cite{helson10} for the precise definition of this weak product space, but recall that
\[\|f\|_\mathscr{K} = \inf_{J<\infty}\left\{\sum_{j=1}^J \|g_j\|_{\mathscr{H}^2}\|h_j\|_{\mathscr{H}^2}\,:\, f(s) = \sum_{j=1}^J g_j(s)h_j(s)\right\}.\]
By the Cauchy--Schwarz inequality, we have $\|f\|_{\mathscr{H}^1}\leq\|f\|_{\mathscr{K}}$, and hence $\mathscr{K}\subseteq\mathscr{H}^1$. Nehari's Theorem \cite{nehari} states that we have $H^1(\mathbb{D})=H^2(\mathbb{D})\odot H^2(\mathbb{D})$. It was shown by Ortega-Cerd\`{a} and Seip in \cite{ortega2012lower} that the corresponding result is false for $\mathscr{K}$ and $\mathscr{H}^1$. Hence the following result \emph{does not} imply that $\mathscr{C}_\Phi$ maps $\mathscr{H}^1$ to $\mathscr{H}^1$.

\begin{thm}
	Let $\Phi \in \mathscr{G}$ with $c_0=0$. Then $\mathscr{C}_\Phi \colon \mathscr{K} \to \mathscr{H}^1$ and $\mathscr{C}_\Phi \colon \mathscr{H}^1 \to \mathscr{A}^2$.
	\begin{proof}
		For the first statement, the argument of \cite{gordonhedenmalm} applies line for line, provided we can show that $\mathscr{K}$ is locally embedded in $H^1(\mathbb{C}_{1/2})$: For $\tau \in \mathbb{R}$ and a Dirichlet polynomial $f$, we have
		\begin{equation} \label{eq:helsonembed} 
			\int_{\tau}^{\tau+1} \left|f\left(\frac{1}{2}+it\right)\right|\,dt \leq C\|f\|_\mathscr{K}. 
		\end{equation}
		The embedding \eqref{eq:helsonembed} follows immediately from the local embedding of $\mathscr{H}^2$ into $H^2(\mathbb{C}_{1/2})$ \cite{HLS} and the Cauchy--Schwarz inequality. For the second statement we shall use Helson's inequality \cite{helson06}: For $f$ of the form \eqref{eq:diriseri}, we have 
		\begin{equation} \label{eq:helson}
			\left(\sum_{n=1}^\infty \frac{|a_n|^2}{d(n)}\right)^\frac{1}{2} \leq \|f\|_{\mathscr{H}^1}.
		\end{equation}
		Clearly, since $\mathscr{A}^2=\mathscr{D}_1$, this means that $\|f\|_{\mathscr{A}^2}\leq\|f\|_{\mathscr{H}^1}$, and hence the second statement follows from Theorem~\ref{thm:compo}.
	\end{proof}
\end{thm}

\begin{rem}
	The space $\mathscr{K}$ is interesting in its own right: Its dual is isometrically isomorphic to the space of all bounded multiplicative Hankel forms on $\mathscr{H}^2$. We refer to \cite{helsonbook}.
\end{rem}

Helson's inequality \eqref{eq:helson} is the $p=1$ case from the family of inequalities 
\[\|f\|_{\mathscr{A}^{2p}}\leq \|f\|_{\mathscr{H}^p}.\] 
\noindent These inequalities are obtained in \cite{BL2013} by an iterative process similar to the one used in \cite{helson06} from the following result: For $1\leq p < \infty$ and $F\in H(\mathbb{D})$, 
\begin{equation} \label{eq:carleman}
	\|F\|_{A^{2p}}=\left(\int_{\mathbb{D}} |F(w)|^{2p}\,dm(w)\right)^\frac{1}{2p} \leq \left(\int_0^{2\pi} \left|F\left(e^{i\theta}\right)\right|^p\, \frac{d\theta}{2\pi}\right)^\frac{1}{p} = \|f\|_{H^p(\mathbb{D})}.
\end{equation}

The inequality \eqref{eq:carleman} has been rediscovered several times, see \cite{vukotic} and \cite{mateljevic}, and dates at least back to Carleman \cite{carleman}. It is essential that the measures are normalized and the constant is $1$ for the iterative procedure to work.

Now, we define $\mathscr{H}^{2p}\odot \mathscr{H}^{2p}$ for $p=1,\,3,\,5,\,7\,\ldots$ in a similar manner as above. Clearly, the space $\mathscr{H}^{2p}\odot\mathscr{H}^{2p}$ is contained in $\mathscr{H}^p$ and $\mathscr{H}^{2p}\odot\mathscr{H}^{2p}$ is moreover locally embedded in $H^p(\mathbb{C}_{1/2})$. Arguing as above, but replacing Theorem~\ref{thm:compo} with Theorem~\ref{thm:Apcomp} for $p=2k$ in the second statement, we obtain the following result:

\begin{cor}
	Let $\Phi \in \mathscr{G}$ with $c_0=0$ and suppose that $p$ is an odd integer. Then $\mathscr{C}_\Phi \colon \mathscr{H}^{2p}\odot\mathscr{H}^
	{2p} \to \mathscr{H}^p$ and $\mathscr{C}_\Phi \colon \mathscr{H}^p \to \mathscr{A}^{2p}$.
\end{cor}

\subsection{} This paper has been devoted to the study of \emph{Bohr--Bergman} spaces of Dirichlet series, that is spaces defined by a Bergman--type norm through the Bohr lift to the polydisc $\mathbb{D}^\infty$. 

As indicated by the fact that a formula similar to \eqref{eq:Hpnorm} is not possible for $\mathbb{D}^\infty$, the Bohr--Bergman spaces are not the only Bergman spaces of Dirichlet series. One could also define limit Bergman-type norms in the half-plane $\mathbb{C}_0$, and take the closure of Dirichlet polynomials with respect to such a norm. 

The main example of such spaces are the Hilbert spaces introduced by McCarthy \cite{mccarthy}, which we for $\alpha>0$ define as
\[\mathcal{A}_\alpha = \left\{f(s) = \sum_{n=1}^\infty a_n n^{-s} \,:\, \|f\|_{\mathcal{A}_\alpha} = \left(\sum_{n=1}^\infty \frac{|a_n|^2}{(1+ \log{n})^\alpha}\right)^\frac{1}{2}<\infty\right\}.\]
By Theorem~1 and Example~1 of \cite{olsenlocal}, we may prove a version of Lemma~\ref{lem:localembed} for $\mathcal{A}_\alpha$, but in this case $\beta=\alpha$. It is clear that the integer $n$ has at most $\log{n}/\log{2}$ prime factors, and hence
\[\Omega(n)+1 \leq \frac{\log{n}}{\log{2}}+1 \leq 2(\log{n}+1).\] 
Combining these observations with \eqref{eq:Omegapull} we obtain the following result: If $\Phi \in \mathscr{G}$ with characteristic $c_0=0$, then $\mathscr{C}_\Phi$ maps $\mathcal{A}_\alpha$ to $\mathcal{A}_\alpha$, for every $\alpha>0$. Composition operators on $\mathcal{A}_\alpha$ have already been studied in \cite{bailleul}. Combining this result with Th{\'e}or{\`e}me 9.1 of \cite{bailleul} (or Theorem~1 in \cite{B2014}) yields a complete description of composition operators on $\mathcal{A}_\alpha$:
\begin{thm} \label{thm:mccomp}
	Let $\alpha>0$. A function $\Phi \colon \mathbb{C}_{1/2} \to \mathbb{C}_{1/2}$ defines a composition operator $\mathscr{C}_\Phi \colon \mathcal{A}_\alpha \to \mathcal{A}_\alpha$ if and only if $\Phi\in\mathscr{G}$. 
\end{thm}
The fact that our ``polydisc point of view'' argument applies to the ``half-plane point of view'' could likely be viewed as a coincidence. In fact, by the irregularity of $d(n)$ it is clear that $\mathcal{A}_\alpha \not\subseteq \mathscr{D}_\beta$ and $\mathscr{D}_\beta \not\subseteq \mathcal{A}_\alpha$ for every $\alpha,\beta>0$.

% ACKNOWLEDGEMENTS
\section*{Acknowledgements}
The authors would like to extend their gratitude to H. Queff{\'e}lec for making them aware of their independent investigations into this topic, thereby inducing a collaboration resulting in the present work. The authors would also like to thank the anonymous referee for providing several pertinent suggestions, which improved the quality of the paper.

% REFERENCES
\bibliographystyle{amsplain} 
\bibliography{comp} 
\end{document}